\documentclass{gtmon_a}
\pdfoutput=1

\proceedingstitle{The interaction of finite-type and Gromov--Witten
invariants (BIRS 2003)}
\conferencestart{15 November 2003}
\conferenceend{20 November 2003}
\conferencename{The interaction of finite-type and Gromov--Witten
invariants}
\conferencelocation{Banff International Research Station, Banff, Alberta,
Canada}

\editor{David Auckly}
\givenname{David}
\surname{Auckly}

\editor{Jim Bryan}
\givenname{Jim}
\surname{Bryan}

\title{Formulae of one-partition and two-partition Hodge integrals}

\author{Chiu-Chu Melissa Liu}
\givenname{Chiu-Chu Melissa}
\surname{Liu}
\address{Department of Mathematics\\
Harvard University\\\newline
1 Oxford Street\\
Cambridge\\MA 02138\\USA}
\email{ccliu@math.harvard.edu}
\urladdr{}

\volumenumber{8}
\issuenumber{}
\publicationyear{2006}
\papernumber{5}
\startpage{105}
\endpage{128}

\doi{}
\MR{}
\Zbl{}


\received{1 January 2005}
\revised{26 January 2005}
\accepted{1 February 2005}
\published{22 April 2006}
\publishedonline{22 April 2006}
\proposed{}
\seconded{}
\corresponding{}
\editor{}
\version{}

\arxivreference{math.AG/0502430}

\begin{asciiabstract}
Prompted by the duality between open string theory on noncompact Calabi-Yau
threefolds and Chern-Simons theory on three-manifolds, M Marino and
C Vafa conjectured a formula of one-partition Hodge integrals in term of
invariants of the unknot. Many Hodge integral identities, including the
lambda_g conjecture and the ELSV formula, can be obtained by taking
limits of the Marino-Vafa formula.

Motivated by the Marino-Vafa formula and  formula of  Gromov-Witten
invariants of local toric Calabi-Yau threefolds predicted by physicists,
J Zhou conjectured a formula of two-partition Hodge integrals in terms
of invariants of the Hopf link and used it to justify the physicists'
predictions.

In this expository article, we describe proofs and applications of
these two formulae of Hodge integrals based on joint works of K Liu,
J Zhou and the author.  This is an expansion of the author's talk of
the same title at the BIRS workshop `The Interaction of Finite
Type and Gromov-Witten Invariants', November 15th--20th 2003.
\end{asciiabstract}

\AtBeginDocument{\let\bar\wbar\let\tilde\wtilde\let\hat\what}

\newcommand{\cM}{\mathcal{M}}
\newcommand{\cO}{\mathcal{O}}
\newcommand{\cW}{\mathcal{W}}

\newcommand{\bC}{\mathbb{C}}
\newcommand{\bE}{\mathbb{E}}
\newcommand{\bL}{\mathbb{L}}
\newcommand{\bP}{\mathbb{P}}
\newcommand{\bQ}{\mathbb{Q}}
\newcommand{\bZ}{\mathbb{Z}}

\newcommand{\Aut}{\mathrm{Aut}}
\newcommand{\Br }{\mathrm{Br} }
\newcommand{\vir}{\mathrm{vir}}
\newcommand{\Sym}{\mathrm{Sym}}

\newcommand{\lam}{\lambda}
\newcommand{\si}{\sigma}

\newcommand{\tf}{\widetilde{f}}

\newcommand{\Mbar}{{}\mskip4mu\overline{\mskip-4mu\mathcal{M}\mskip-1mu}\mskip1mu}
\newcommand{\pa}{\partial}
\newcommand{\bu}{\bullet}

\makeatletter
\def\cnewtheorem#1[#2]#3{\newtheorem{#1}{#3}[section]
\expandafter\let\csname c@#1\endcsname\c@theo}

\newtheorem{theo}{Theorem}
\cnewtheorem{prop}[theo]{Proposition}
\cnewtheorem{lemm}[theo]{Lemma}
\cnewtheorem{rema}[theo]{Remark}
\cnewtheorem{defi}[theo]{Definition}
\cnewtheorem{coro}[theo]{Corollary}
\cnewtheorem{exam}[theo]{Example}
\cnewtheorem{conj}[theo]{Conjecture}
\cnewtheorem{fact}[theo]{Fact}

\theoremstyle{remark}

\makeatother  

\begin{document}

\begin{abstract} 
Prompted by the duality between open string theory on noncompact Calabi--Yau
threefolds and Chern--Simons theory on three-manifolds, M Mari\~{n}o and
C Vafa conjectured a formula of one-partition Hodge integrals in term of
invariants of the unknot. Many Hodge integral identities, including the
$\lambda_g$ conjecture and the ELSV formula, can be obtained by taking
limits of the Mari\~{n}o--Vafa formula.

Motivated by the Mari\~{n}o--Vafa formula and  formula of  Gromov--Witten
invariants of local toric Calabi--Yau threefolds predicted by physicists,
J Zhou conjectured a formula of two-partition Hodge integrals in terms
of invariants of the Hopf link and used it to justify the physicists'
predictions.

In this expository article, we describe proofs and applications of
these two formulae of Hodge integrals based on joint works of K Liu,
J Zhou and the author.  This is an expansion of the author's talk of
the same title at the BIRS workshop \emph{The Interaction of Finite
Type and Gromov--Witten Invariants}, November 15th--20th 2003.
\end{abstract}

\maketitle

\section{Introduction}
\label{sec-introduction}

In \cite{WiTwo}, Witten related topological string theory on the cotangent
bundle $T^*M$ of a three manifold $M$ to the Chern--Simons gauge
theory on $M$.  Gopakumar and Vafa \cite{Go-Va} related the topological
string theory on the deformed conifold (the cotangent bundle $T^*S^3$
of the 3--sphere) to that on the resolved conifold (the total space of
$\cO_{\bP^1}(-1)\oplus \cO_{\bP^1}(-1)\to \bP^1$). These works lead to
a duality between the topological string theory on the resolved conifold
and the Chern--Simons theory on $S^3$. A mathematical consequence of this
duality is a surprising relationship between Gromov--Witten invariants,
which arise in the topological string theory, and knot invariants,
which arise in the Chern--Simons theory.  Ooguri and Vafa \cite{Oo-Va}
proposed that Chern--Simons knot invariants of a knot $K$ in $S^3$ should
be related to open Gromov--Witten invariants of $(X,L_K)$, where $X$
is the resolved conifold and $L_K$ is a Lagrangian submanifold of $X$
canonically associated to~$K$.

Ooguri and Vafa's proposal for a general knot is not precise enough to
be a mathematical conjecture for the following reasons. First of all,
Ooguri and Vafa did not give an explicit description of $L_K$ for a
general knot $K$; later Taubes carried out a construction  for a general
knot \cite{Ta} so we now have a candidate.  Secondly, given a Lagrangian
submanifold $L$ of $X$, it is not known how to define open Gromov--Witten
invariants of $(X,L)$ in general, though these invariants are expected
to count holomorphic curves in $X$ with boundary in $L$.

When $K$ is the unknot, Ooguri and Vafa's proposal can be made more
precise.  In this case, Ooguri and Vafa described $L_K$ explicitly
and conjectured a formula of open Gromov--Witten invariants of
$(X,L_K)$. Although these invariants were not defined, Katz--Liu
\cite{Ka-Li} and Li--Song \cite{Li-So} carried out heuristic localization
calculations which agreed with the Ooguri--Vafa formula. One-partition
Hodge integrals arise in Katz and Liu's calculations with varying
torus weight. Mari\~{n}o and Vafa identified the weight dependence of
one-partition Hodge integrals with the framing dependence of invariants
of the unknot and conjectured a formula of one-partition Hodge integrals
\cite{Ma-Va}. This formula is a precise mathematical statement, and has
strong consequences on Hodge integrals.

Ooguri and Vafa's proposal can be generalized to links, 
Labastida--Mari{\~n}o--Vafa \cite{La-Ma-Va}. In
particular, the Hopf link corresponds to an explicit Lagrangian
submanifold with two connected components, and the relevant  open
Gromov--Witten invariants can be transformed to two-partition Hodge
integrals by heuristic localization  calculations \cite{Di-Fl-Gr}.  This
leads to a formula of two-partition Hodge integrals. Zhou conjectured this
formula in \cite{ZhTwo} and used it to verify the formula of Gromov--Witten
invariants of any local toric Calabi--Yau threefold proposed by physicists,
Aganagic--Mari{\~n}o--Vafa \cite{Ag-Ma-Va},
Iqbal \cite{Iq}, Zhou \cite{ZhThree}.

The discussion above is a brief description of the origin of the
formulae of one-partition and two-partition Hodge integrals. See
Mari{\~n}o \cite{MaOne,MaTwo} for much more complete survey on the
duality between topological string theory on noncompact Calabi--Yau
threefolds and Chern--Simons theory on three manifolds.

In this paper, we will describe proofs and applications of these two
formulae of Hodge integrals based on Liu--Liu--Zhou
\cite{LLZOne,LLZTwo,LLZThree,LLZFour}.  We now give an overview of the
rest of this paper. In \fullref{sec-I}, we define one-partition Hodge
integrals and state the formula of such integrals conjectured by
Mari\~{n}o and Vafa.  In \fullref{sec-applicationI}, we describe how
to extract the values of all $\lam_g$--integrals and the ELSV formula
from the Mari\~{n}o--Vafa formula, following \cite{LLZThree}. In
\fullref{sec-three}, we describe three approaches to the
Mari\~{n}o--Vafa formula: the first approach is based on cut-and-join
equations and is used in the proof of the Mari\~{n}o--Vafa formula
given in \cite{LLZOne,LLZTwo}; the second approach is based on
convolution equations and can be generalized to prove the formula of
two-partition Hodge integrals \cite{LLZFour}; the third approach is
based on bilinear localization equations and is used in Okounkov and
Pandharipande's proof of the Mari\~{n}o--Vafa formula
\cite{Ok-PaThree}.  In \fullref{sec-KGproof}, we present a proof of
the convolution equation for one-partition Hodge integrals.  In
\fullref{sec-II}, we generalize the previous discussion to the
two-partition case.

Finally, this paper is based on the author's talk at the BIRS workshop
\emph{The Interaction of Finite Type and Gromov--Witten Invariants}
and does not cover developments since the workshop took place.  The
two-partition Hodge integrals and their applications can be better
understood in terms of the topological vertex
Aganagic--Klemm--Mari{\~n}o--Vafa \cite{AKMV}.  The topological vertex
is related to three-partition Hodge integrals Diaconescu--Florea
\cite{Di-Fl}, and a mathematical theory is developed in
Li--Liu--Liu--Zhou \cite{LLLZ}.  Another exciting development is the
correspondence between Gromov--Witten invariants and Donaldson--Thomas
invariants Maulik--Nekrasov--Okounkov--Pandharipande
\cite{MNOPOne,MNOPTwo}.

\subsection*{Acknowledgments} It is a great pleasure to thank my
collaborators K Liu, J Zhou, and later J Li. I understand the mathematical
aspect of the duality between Gromov--Witten theory and Chern--Simons
theory much better than I did before I started to collaborate with
them. I also wish to thank D Auckly and J Bryan, the organizers of the
BIRS workshop \emph{The Interaction of Finite Type and Gromov--Witten
Invariants}, for inviting me to give lectures and contribute to the
proceedings.  The workshop directly benefited my research. Finally,
I wish to thank A Okounkov and R Pandharipande for explaining their
works that are mentioned in this paper.

\section{The Mari\~{n}o--Vafa formula of one-partition Hodge integrals}
\label{sec-I}
In this section, we state the formula of one-partition Hodge integrals
conjectured by Mari\~{n}o and Vafa. We first recall the definitions
of partitions and Hodge integrals.

\subsection{Partitions}
\label{sec-partition}
Recall that  a \emph{partition} $\mu$ of a nonnegative integer $d$
is a sequence of positive numbers
$$\mu=(\mu_1\geq\cdots\geq\mu_h >0)$$
such that $\mu_1+\cdots+\mu_h=d$. We call $\ell(\mu)=h$ the \emph{length}
of $\mu$ and $|\mu|=d$ the \emph{size} of $\mu$. The automorphism
group $\Aut(\mu)$ permutes $\mu_i$ and $\mu_j$ if $\mu_i=\mu_j$. For
example, $\Aut(5,5,4,1,1,1)\cong S_2\times S_3$, where $S_n$ denotes
the permutation group of $n$ elements.  In particular, $\Aut(\mu)$ is
trivial if and only if all the components $\mu_1,\ldots,\mu_h$ of $\mu$
are distinct.

\subsection{Hodge integrals}
\label{sec-hodge}
Let $\Mbar_{g,h}$ be the Deligne--Mumford compactification of moduli space of complex algebraic
curves of genus $g$ with $h$ marked points. A point in $\Mbar_{g,h}$ is represented by
$(C,x_1,\ldots,x_h)$, where $C$ is a complex algebraic curve of arithmetic genus $g$ with at most
nodal singularities, $x_1,\ldots, x_h$ are distinct smooth points on $C$, and $(C,x_1,\ldots,x_h)$ is
\emph{stable} in the sense that its automorphism group is finite.
  
The Hodge bundle $\bE$ is a rank--$g$ vector bundle over
$\Mbar_{g,h}$ whose fiber over the moduli point
$$[(C,x_1,\ldots,x_h)]\in \Mbar_{g,h}$$
is $H^0(C,\omega_C)$. When $C$ is smooth, it can be viewed as a compact Riemann surface and
$H^0(C,\omega_C)$ is the space of holomorphic one forms on $C$.
The $\lam$--classes are defined by
$$\lam_j=c_j(\bE)\in H^{2i}\bigl(\Mbar_{g,h};\bQ\bigr).$$
The cotangent line $T_{x_i}^* C$ of $C$ at the $i$th marked point $x_i$
gives rise to a line bundle $\bL_i$ over $\Mbar_{g,h}$.
The $\psi$--classes are defined by
$$\psi_i=c_1(\bL_i)\in H^2\bigl(\Mbar_{g,h};\bQ\bigr).$$
The $\lam$--classes and $\psi$--classes lie in $H^*(\Mbar_{g,h};\bQ)$
instead of $H^*(\Mbar_{g,h};\bZ)$ because
$\bE$ and $\psi$ are orbibundles on the compact orbifold $\Mbar_{g,h}$.

\emph{Hodge integrals} are intersection numbers of $\lambda$--classes and
$\psi$--classes:
$$\int_{\Mbar_{g,h}}\psi_1^{j_1}\cdots \psi_h^{j_h} 
  \lam_1^{k_1}  \cdots \lam_g^{k_g} \in \bQ$$
The $\psi$--integrals (also known as \emph{descendent integrals})
$$\int_{\Mbar_{g,h}}\psi_1^{j_1}\cdots \psi_h^{j_h}$$
can be computed recursively by Witten's conjecture \cite{WiOne}
proved by Kontsevich \cite{KoOne}. Okounkov and Pandharipande gave
another proof of Witten's conjecture in \cite{Ok-PaOne}, and M Mirzakhani
recently gave a third proof \cite{Mi}.

Using Mumford's Grothendieck--Riemann--Roch calculations \cite{Mu}, Faber
and Pandharipande showed \cite{Fa-PaOne} that general Hodge integrals can
be uniquely reconstructed from descendent integrals.  See \cite[Section
4]{Gi} for Givental's reformulation of this result.

\subsection{One-partition Hodge integrals}
\label{sec-hodgeI}
Given a triple $(g,\mu,\tau)$, where $g$ is a nonnegative
integer, $\mu$ is a partition, and $\tau\in\bZ$, define
a \emph{one-partition Hodge integral} $G_{g,\mu}(\tau)$ as
\begin{align*}
\frac{-\sqrt{-1}^{|\mu|+\ell(\mu)}(\tau(\tau{+}1))^{\ell(\mu)-1}}{|\Aut(\mu)|}
&\Biggl(\,\prod_{i=1}^{\ell(\mu)}\frac{\prod_{a=1}^{\mu_i-1}(\mu_i\tau{+}a)}{(\mu_i{-}1)!}\Biggr)\\
&\hspace{2cm}\int_{\Mbar_{g,\ell(\mu)}}\!\!\!\!\!\!
\frac{\Lambda_g^\vee(1)\Lambda_g^\vee(-\tau{-}1)\Lambda_g^\vee(\tau)}
{\prod_{i=1}^{\ell(\mu)}(1{-}\mu_i\psi_i)}
\end{align*}
$$
\Lambda_g^\vee(u)=u^g-\lam_1 u^{g-1}+\ldots+(-1)^g\lam_g.\leqno{\hbox{where}}
$$
One-partition Hodge integrals can be simplified in special cases.
If $g=0$, then $\Lambda^\vee_0(u)=1$ and
\begin{align*}
\int_{\Mbar_{0,h}}\frac{1}{\prod_{i=1}^h(1-\mu_i\psi_i)}
&=\sum_{k_1+\cdots+k_h=h-3}\mu_1^{k_1}\cdots\mu_h^{k_h}
\int_{\Mbar_{0,h}}\psi_1^{k_1}\cdots\psi_h^{k_h}\\
&=\sum_{k_1+\cdots+k_h=h-3}\mu_1^{k_1}\cdots\mu_h^{k_h}
    \frac{h!}{k_1!\cdots k_h!}\\
&=|\mu|^{h-3}.
\end{align*}
If $\tau=0$, then $G_{g,\mu}(0)=0$ if $\ell(\mu)>1$, and
$$G_{g,(d)}(0)=-\sqrt{-1}^{d+1}\int_{\Mbar_{g,1}}\frac{\lambda_g}{1-d\psi}
=-\sqrt{-1}^{d+1}d^{2g-2}b_g,$$
$$b_g=\left\{\begin{array}{ll}1,& g=0,\\
\int_{\Mbar_{g,1}}\lambda_g\psi^{2g-2}, & g>0.\end{array}\right.\leqno{\hbox{where}}$$

\subsection{Formula of one-partition Hodge integrals}
\label{sec-formulaI}
To state the formula of one-partition Hodge integrals, we need
to introduce some generating functions.

We first define generating functions of one-partition Hodge integrals.
Introduce variables $\lambda$ and $p=(p_1,p_2,\ldots)$.
Given a partition $\mu$, let
$$
p_\mu=p_{\mu_1}\ldots p_{\mu_{\ell(\mu)} }.
$$
Define generating functions
\begin{align*}
G_\mu(\lam;\tau)&=
\sum_{g=0}^\infty\lam^{2g-2+ \ell(\mu)}G_{g,\mu}(\tau),\\
G(\lam;\tau;p)&=\sum_{\mu\neq \emptyset}G_\mu(\lam;\tau)p_\mu,\\
G^\bu(\lam;\tau;p)&= \exp\left(G(\lam;\tau;p)\right)
= \sum_\mu G^\bu_\mu (\lam;\tau)p_\mu
= 1+\sum_{\mu\neq\emptyset} G^\bu_\mu(\lam;\tau)p_\mu.
\end{align*}
Here $\emptyset$ is the empty partition, the unique partition such that
$|\emptyset|=\ell(\emptyset)=0$, and 
$G^\bu_\mu(\lam;\tau)$ is the disconnected version of $G_\mu(\lam;\tau)$.

We next define generating functions of symmetric group representations.
Introduce
\begin{equation}\label{eqn-WI}
\cW_{\mu}(q)= q^{\kappa_\mu/4} \prod_{1 \leq i < j \leq \ell(\mu)}\frac{[\mu_i-\mu_j+j-i]}{[j-i]}
\prod_{i=1}^{\ell(\mu)}\frac{1}{\prod_{v=1}^{\mu_i}[v-i+\ell(\mu)] }
\end{equation}
where 
$$
\kappa_{\mu} = |\mu| + \sum_i (\mu_i^2 - 2i\mu_i),\qquad
[m]=q^{m/2}-q^{-m/2},\qquad
q=e^{\sqrt{-1}\lam}.
$$ 
The expression $\cW_{\mu}(q)$ is related to the HOMFLY  polynomial of the unknot.  
Let  $\chi_\mu$ denote the character of the irreducible representation of
symmetric group $S_{|\mu|}$ indexed by $\mu$,
and let $C_\mu$ denote the conjugacy class of $S_{|\mu|}$ indexed by
$\mu$.  Define
\begin{equation}\label{eqn-RV}
R^\bu_\mu(\lam; \tau) = \sum_{|\nu|=|\mu|}
\frac{\chi_{\nu}(C_\mu)}{z_{\mu}}e^{\sqrt{-1}\tau\kappa_{\nu}\lam/2}
\sqrt{-1}^{|\mu|}\cW_{\nu}(q),
\end{equation}
where $ z_\mu = |\Aut(\mu)|\mu_1\cdots\mu_{\ell(\mu)}$.
Finally, define
$$
R^\bu(\lam;\tau;p)=\sum_\mu R^\bu_\mu(\lam;\tau)p_\mu
$$
and its connected version
$$
R(\lam;\tau;p)=\log R^\bu(\lam;\tau;p)=\sum_{\mu\neq\emptyset} R_\mu(\lam;\tau)p_\mu.
$$
\begin{conj}[{Mari\~{n}o--Vafa \cite{Ma-Va}}]
\begin{equation}\label{eqn-MV}
G(\lam;\tau;p) = R(\lam;\tau;p).
\end{equation}
\end{conj}

The Mari\~{n}o--Vafa formula \eqref{eqn-MV} provides a highly
nontrivial link between geometry (Hodge integrals) and combinatorics
(representations of symmetric groups). Note that for each fixed partition
$\mu$, the Mari\~{n}o--Vafa formula gives a closed and finite formula
of $G_\mu(\lam;\tau)$, a generating function of all genera.

\section{Applications of the Mari\~{n}o--Vafa formula}
\label{sec-applicationI}
Many Hodge integral identities can be obtained by taking limits
of the Mari\~{n}o--Vafa formula. We illustrate this by some examples,  
following \cite{LLZThree}.

We have
$$
G_{g,\mu}(\tau)=\sum_{k=\ell(\mu)-1}^{2g-2+|\mu|+\ell(\mu)} G_{g,\mu}^k \tau^k,
$$
where
\begin{align*}
G_{g,\mu}^{\ell(\mu)-1}&=
\frac{-\sqrt{-1}^{|\mu|+\ell(\mu)}}{|\Aut(\mu)|}
\frac{\lambda_g}{\prod_{i=1}^{\ell(\mu)}(1-\mu_i\psi_i) },\\
G_{g,\mu}^{2g-2+|\mu|+\ell(\mu)}&=
\frac{-\sqrt{-1}^{|\mu|+\ell(\mu)}}{|\Aut(\mu)|}
\Biggl(\,\prod_{i=1}^{\ell(\mu)}\frac{\mu_i^{\mu_i}}{\mu_i!}\,\Biggr)
\int_{\Mbar_{g,\ell(\mu)}}
\frac{\Lambda_g^\vee(1)}{\prod_{i=1}^{\ell(\mu)}(1-\mu_i\psi_i) }.
\end{align*}
\subsection{$\lam_g$--integrals}
\label{sec-lambda}
Extracting the part corresponding to $G^{\ell(\mu)-1}_{g,\mu}$ from
$R(\lam;\tau;p)$, we obtain
\begin{equation}\label{eqn-pre}
\sum_{g=0}^\infty \lam^{2g}
\int_{\Mbar_{g,n}}\frac{\lam_g}{\prod_{i=1}^n(1-\mu_i\psi_i)}
=|\mu|^{n-3}\frac{|\mu|\lam/2}{\sin(|\mu|\lam/2)}.
\end{equation}
The identity \eqref{eqn-pre} is true for any partition of length $n$, so we may view
it as an identity in $\bQ[\mu_1,\ldots,\mu_n][[\lam]]$. This gives
us the values of all $\lambda_g$--integrals:
\begin{equation}\label{eqn-lambda}
\int_{\Mbar_{g,n}}\psi_1^{k_1}\cdots \psi_n^{k_n}\lambda_g
= {2g+n-3 \choose k_1,\cdots,k_n }\frac{2^{2g-1}-1}{2^{2g-1}}\frac{|B_{2g}|}{(2g)!}
\end{equation}
where $B_{2g}$ are Bernoulli numbers.

Equation \eqref{eqn-lambda} also follows from the following two identities:
\begin{fact}[A formula for $b_g$]\label{thm-bg}
\begin{equation}\label{eqn-gone}
\sum_{g=0}^\infty b_g t^{2g}=\frac{t/2}{\sin(t/2)}
\end{equation}
\end{fact}
\begin{fact}[The $\lam_g$ conjecture]\label{thm-lamg}
\begin{equation}\label{eqn-gn}
\int_{\Mbar_{g,n}}\psi_1^{k_1}\cdots \psi_n^{k_n}\lam_g
= {2g+n-3 \choose k_1,\cdots,k_n }b_g
\end{equation}
\end{fact}
Recall that $b_g$ are defined in \fullref{sec-hodgeI}.  Formula \eqref{eqn-gone} gives 
values of $b_g$ which are $\lambda_g$--integrals on $\Mbar_{g,1}$, while
the $\lam_g$ conjecture tells us how general $\lam_g$--integrals are
determined by $b_g$.
The formula \eqref{eqn-gone} was proved by Faber and Pandharipande in \cite{Fa-PaOne}.
The $\lam_g$ conjecture was found by Getzler and Pandharipande
\cite{Ge-Pa} as a consequence of the
degree--$0$ Virasoro conjecture of $\bP^1$ and 
was first proved by Faber and Pandharipande \cite{Fa-PaTwo}. 
Later the Virasoro conjecture was proved for projective spaces \cite{Gi}
and curves \cite{Ok-PaFour}; both cases include $\bP^1$ as a special case.

\subsection{ELSV formula}
\label{sec-ELSV}
The part corresponding to $G_{g,\mu}^{2g-2+|\mu|+\ell(\mu)}$
in $R(\lam;\tau;p)$ reduces to the Burnside formula of Hurwitz numbers.
Recall that the Hurwitz number  $H_{g,\mu}$ is the weighted count of
genus--$g$, degree--$|\mu|$ ramified covers of $\bP^1$ with fixed
ramification type $\mu$ over a point $q^1\in\bP^1$.
We obtain the ELSV formula:
\begin{equation}\label{eqn-ELSV}
\frac{1}{|\Aut(\mu)|}
\Biggl(\,\prod_{i=1}^{\ell(\mu)}\frac{\mu_i^{\mu_i}}{\mu_i!}\,\Biggr)
\int_{\Mbar_{g,\ell(\mu)}}
\frac{\Lambda_g^\vee(1)}{\prod_{i=1}^{\ell(\mu)}(1-\mu_i\psi_i) }
=\frac{H_{g,\mu}}{(2g-2+|\mu|+\ell(\mu))!}.
\end{equation}
This identity was first proved by Ekedahl, Lando, Shapiro, and Vainshtein
\cite{ELSVOne,ELSVTwo}.  In \cite{Gr-VaOne}, T Graber and R Vakil gave
a proof by virtual localization on moduli spaces $\Mbar_g(\bP^1,d)$
of stable maps to $\bP^1$, and described a simplified proof by virtual
localization on moduli spaces $\Mbar_{g,0}(\bP^1,\mu)$ of relative
stable maps to $(\bP^1,q^1)$ (see \fullref{sec-moduli} for precise
definition).  Actually, virtual localization on $\Mbar_{g,0}(\bP^1,\mu)$
with suitable choices of weights gives both the ELSV formula and the
cut-and-join equation of Hurwitz numbers \cite[Section 7]{LLZTwo}. The
latter was first proved using combinatorics in \cite{Go-Ja-Va}, and later
using the symplectic sum formula in \cite{Li-Zh-Zh} and \cite[Section
15.2]{Io-PaTwo}.

\subsection{Other identities}
\label{sec-identity}
The following identities, proved in \cite{Fa-PaOne},
are also consequences the Mari\~{n}o--Vafa formula:
\begin{eqnarray*}
\int_{\Mbar_g}\lambda_{g-2}\lambda_{g-1}\lambda_g&=&
\frac{1}{2(2g-2)!}\frac{|B_{2g-2}|}{2g-2}\frac{|B_{2g}|}{2g}\\
\int_{\Mbar_{g,1}}\frac{\lambda_{g-1}}{1-\psi_1}
&=&b_g\sum_{i=1}^{2g-1}\frac{1}{i}-\frac{1}{2}
\sum_{\substack{g_1+g_2=g\\ g_1,g_2>0}}
\frac{(2g_1-1)!(2g_2-1)!}{(2g-1)!}b_{g_1}b_{g_2}
\end{eqnarray*}
See \cite{LLZThree} for details.

\section{Three approaches to the Mari\~{n}o--Vafa formula}
\label{sec-three}
\subsection{The cut-and-join equation and functorial localization}
\label{sec-functorial}

In this subsection, we described the strategy of the proof of 
the Mari\~{n}o--Vafa formula given in \cite{LLZTwo}.
We have seen in \fullref{sec-hodgeI} that the left hand side of
the Mar\~{n}o--Vafa formula can be greatly simplified at $\tau=0$:
\begin{equation}\label{eqn-Ginitial}
G(\lambda;0;p)=-\sum_{d=1}^\infty\frac{\sqrt{-1}^{d+1}p_d}{\lambda d^2}
\sum_{g=0}^\infty b_g(\lambda d)^{2g}
\end{equation}
It turns out that the right hand side of the Mari\~{n}o--Vafa formula can
also be greatly simplified at $\tau=0$ (see \cite{ZhOne} or \cite[Section
2]{LLZTwo} for details):
\begin{equation}\label{eqn-Rinitial}
R(\lambda;0;p)=-\sum_{d=1}^\infty\frac{\sqrt{-1}^{d+1}p_d}
{2d\sin(\lambda d/2)}
\end{equation}
Expressions \eqref{eqn-Ginitial} and \eqref{eqn-Rinitial} are equal by
\fullref{thm-bg}, so the Mari\~{n}o--Vafa
formula \eqref{eqn-MV} holds at $\tau=0$:
\begin{equation}\label{eqn-MVinitial}
G(\lambda;0;p)=R(\lambda;0;p)
\end{equation}
Note that both sides of the Mari\~{n}o--Vafa formula are valid for
$\tau\in\bC$.  From \eqref{eqn-RV} and the orthogonality
of characters
\begin{equation}\label{eqn-orthogonal}
\sum_{\si}\frac{\chi_\mu(C_\si)\chi_\nu(C_\si)}{z_\si}=\delta_{\mu\nu}
\end{equation}
it follows immediately that
\begin{equation}\label{eqn-Rcon}
R_\mu^\bu (\lam;\tau)=\sum_{|\nu|=|\mu|} 
R^\bu_\nu(\lam;0)z_\nu\Phi^\bu_{\nu,\mu}(\sqrt{-1}\lambda\tau)
\end{equation} 
where
\begin{equation}\label{eqn-burnside}
\Phi^\bu_{\nu,\mu}(\lam)=\sum_{\eta}
\frac{\chi_\eta(C_\nu)}{z_\nu}\frac{\chi_\eta(C_\mu)}{z_\mu} e^{\kappa_\eta\lambda/2}.
\end{equation}
The convolution equation \eqref{eqn-Rcon} is equivalent to the following
\emph{cut-and-join equation}:
\begin{equation}\label{eqn-Rcj}
\frac{\pa R}{\pa \tau}
=\frac{\sqrt{-1}\lambda}{2}\sum_{i,j=1}^\infty
\left((i+j)p_ip_j\frac{\pa R}{\pa p_{i+j} }
+ij p_{i+j}\left( \frac{\pa R}{\pa p_i}
\frac{\pa R}{\pa p_j} +\frac{\pa^2 R}{\pa p_i\pa p_j}\right)\right)
\end{equation}
See \cite{Go-Ja-Va}, \cite{ZhOne}, and \cite[Section 3]{LLZTwo} for details.

The Mari\~{n}o--Vafa formula will follow from the initial values 
\eqref{eqn-MVinitial}, the cut-and-join equation \eqref{eqn-Rcj}
of $R(\lam;\tau;p)$,
and the following cut-and-join equation of $G(\lam;\tau;p)$.
\begin{theo}[{Liu--Liu--Zhou \cite{LLZTwo}}]\label{thm-Gcj}
\begin{equation}\label{eqn-Gcj}
\frac{\pa G}{\pa \tau}
=\frac{\sqrt{-1}\lambda}{2}\sum_{i,j=1}^\infty
\left((i+j)p_ip_j\frac{\pa G}{\pa p_{i+j} }
+ij p_{i+j}\left( \frac{\pa G}{\pa p_i}
\frac{\pa G}{\pa p_j} +\frac{\pa^2 G}{\pa p_i\pa p_j}\right)\right)
\end{equation}
\end{theo}

In \cite{LLZTwo}, \fullref{thm-Gcj} was proved by applying (virtual)
functorial localization \cite{LLY} to the branch morphism
$$\Br\co\Mbar_{g,0}(\bP^1,\mu)\to \Sym^r\bP^1\cong \bP^r,$$
where $\Mbar_{g,0}(\bP^1, \mu)$ is the moduli space of
relative stable maps from a genus--$g$ curve to $\bP^1$
with fixed ramification type $\mu=(\mu_1,\ldots,\mu_h)$ over $q^1\in \bP^1$,
and
$$r=2g-2+|\mu|+\ell(\mu)$$
is the virtual dimension of $\Mbar_{g,0}(\bP^1,\mu)$. The precise definition
of $\Mbar_{g,0}(\bP^1,\mu)$ is given in \fullref{sec-moduli}.
Note that the $\bC^*$--action on $\bP^1$ induce $\bC^*$--actions
on the domain and the target of $\Br$, and $\Br$ is $\bC^*$--equivariant.
This proof was outlined in \cite{LLZOne} and presented in detail in
\cite{LLZTwo}.

\subsection{Convolution equation and double Hurwitz numbers}
\label{sec-dHurwitz}
We now describe a variant of the above approach, which is even more
direct and can be generalized to prove the formula of two-partition
Hodge integrals \cite{LLZFour}. This alternative proof of the
Mari\~{n}o--Vafa formula was discovered by the authors of \cite{LLZTwo}
while working on \cite{LLZFour}.

The Mari\~{n}o--Vafa formula will follow from the initial values \eqref{eqn-MVinitial}, 
the convolution equation \eqref{eqn-Rcon} of $R^\bu_\mu(\lam;\tau)$, 
and the following convolution equation \eqref{eqn-Gcon} of $G^\bu_\mu(\lam;\tau)$:
\begin{theo}[Liu--Liu--Zhou]\label{thm-Gcon}
\begin{equation}\label{eqn-Gcon}
G_\mu^\bullet(\lambda;\tau)
=\sum_{|\nu|=|\mu|} G_\nu^\bullet(\lambda;0)z_\nu
\Phi^\bullet_{\nu,\mu}(\sqrt{-1}\lambda\tau).
\end{equation}
\end{theo}
Recall that the $\lambda_g$ conjecture \eqref{eqn-gn} tells us how
$b_g$ ($\lambda_g$--integrals on $\Mbar_{g,1}$) determine general
$\lambda_g$--integrals.  The convolution equation \eqref{eqn-Gcon}
tells us how the $b_g$ determine all one-partition Hodge integrals (which
contain all $\lambda_g$--integrals and more), so it can be viewed as a
generalization of the $\lambda_g$ conjecture.

By the Burnside formula of Hurwitz numbers, $\Phi^\bu_{\nu,\mu}(\lam)$
is the generating function of disconnected double Hurwitz numbers
$H^\bu_{\chi,\nu,\mu}$:
\begin{equation}\label{eqn-H}
\Phi_{\nu,\mu}^\bu(\lam)=\sum_{\chi}\lam^{-\chi+\ell(\nu)+\ell(\mu)}
\frac{H^\bu_{\chi,\nu,\mu}}{(-\chi+\ell(\nu)+\ell(\mu))!}
\end{equation}
where $|\nu|=|\mu|=d$, and $H^\bu_{\chi,\nu,\mu}$ is the weighted count of
degree--$d$
covers of $\bP^1$ with prescribed ramification types $\nu,\mu$ over two points 
$q^0,q^1\in\bP^1$ by possibly disconnected Riemann surfaces of Euler characteristic $\chi$. 
From \eqref{eqn-burnside} it is clear that 
\begin{equation}\label{eqn-Hsum}
\sum_{|\si|=d}\Phi^\bu_{\nu,\si}(\lam_1)z_\si \Phi^\bu_{\si,\mu}(\lam_2)
=\Phi^\bu_{\nu,\mu}(\lam_1+\lam_2)\quad\text{and}\quad
\Phi^\bu_{\nu,\mu}(0)=\frac{\delta_{\nu,\mu}}{z_\nu}. 
\end{equation}
In \fullref{sec-KGproof}, we will define a generating function
$K^\bu_\mu(\lam)$ of relative Gromov--Witten invariants of $(\bP^1,\infty)$. 
By virtual localization, $K^\bu_\mu(\lam)$ can be expressed in terms of 
one-partition Hodge integrals and double Hurwitz numbers:
\begin{prop}\label{thm-KGI}
\begin{equation}\label{eqn-KGI}
K^\bu_\mu(\lam)=\sum_{|\nu|=|\mu|}G_\nu^\bu(\lambda;\tau)z_\nu
\Phi^\bu_{\nu,\mu}\bigl(-\sqrt{-1}\lambda\tau\bigr).
\end{equation}
\end{prop}
\fullref{thm-KGI} is a special case of \cite[Proposition
7.1]{LLZFour} and, by \eqref{eqn-Hsum}, is equivalent to
\begin{equation}\label{eqn-GKI}
G^\bu_\mu(\lam;\tau)=\sum_{|\nu|=|\mu|}K_\nu^\bu(\lam)z_\nu
\Phi^\bu_{\nu,\mu}\bigl(\sqrt{-1}\lam\tau\bigr),\qquad
G^\bu_\mu(\lam;0)=K^\bu_\mu(\lam).
\end{equation} 
This gives \fullref{thm-Gcon}.

\subsection{Bilinear localization equations}
\label{sec-OP}
This approach is due to Okounkov and Pandharipande \cite{Ok-PaThree}
and was motivated by Faber and Pandharipande's proof of
the $\lam_g$ conjecture \cite{Fa-PaTwo}.

By virtual localization on $\Mbar_{g,n}(\bP^1,d)$,
Okounkov and Pandharipande derived homogeneous bilinear
equations of the form
\begin{equation}\label{eqn-hom}
\sum(\mathrm{linear\ Hodge\ integrals})\cdot
(\mathrm{special\ cubic\ Hodge\ integrals})=0
\end{equation}
and inhomogeneous bilinear equations of the form
\begin{equation}\label{eqn-inhom}
\sum (\mathrm{linear\ Hodge\ integrals})\cdot
(\mathrm{special\ cubic\ Hodge\ integrals})
=(\lam_g\ \mathrm{integrals})
\end{equation}
where ``linear Hodge integrals'' are those in the ELSV formula, 
and ``special cubic Hodge integrals'' are those in the
Mari\~{n}o--Vafa formula. The values of linear Hodge integrals are given
by the ELSV formula \eqref{eqn-ELSV} and the Burnside formula of Hurwitz
numbers; the values of $\lambda_g$--integrals are given by
\eqref{eqn-lambda}. Therefore, \eqref{eqn-hom} and \eqref{eqn-inhom} can be 
viewed as linear equations satisfied by special cubic Hodge integrals.
It was shown in \cite{Ok-PaFour} that there is a unique solution to this system of
linear equations, and  the Mari\~{n}o--Vafa formula gives a solution.

\section[Proof of Proposition \ref{thm-KGI}]
{Proof of \fullref{thm-KGI}}\label{sec-KGproof}

We fix a degree $d$, and a partition $\mu=(\mu_1,\ldots,\mu_h)$ of $d$.

\subsection{Moduli spaces}\label{sec-moduli}
Let $\cM_{g,0}(\bP^1,\mu)$ be the moduli space 
of ramified covers
$$f\co(C,x_1,\ldots, x_h) \to (\bP^1,q^1)$$
of degree--$d$ from a smooth curve $C$ of genus $g$ to $\bP^1$ such that the
ramification type over a distinguished point $q^1\in \bP^1$ is specified
by the partition $\mu$, that is,
$$f^{-1}(q^1)=\mu_1 x_1+\cdots +\mu_h x_h$$ 
as Cartier divisors. 
The moduli space $\cM_{g,0}(\bP^1,\mu)$ is not compact. To compactify it, we
consider the moduli space
$$
\Mbar_{g,0}(\bP^1,\mu)
$$
of relative stable maps to $(\bP^1,q^1)$. 
The moduli spaces of relative stable maps were constructed by
A Li--Ruan \cite{Li-Ru} and Ionel--Parker \cite{Io-PaOne,Io-PaTwo}
in symplectic geometry. Later J Li carried out the construction in
algebraic geometry \cite{LiOne,LiTwo}. We need to use J Li's algebraic 
version. Such  moduli spaces are defined for a general pair $(X,D)$ where
$X$ is a smooth projective variety and $D$ is a smooth divisor. Here
we content ourselves with the definition for this special case.

We first introduce some notation. Let $\bP^1(m)=\bP^1_1\cup\cdots\cup
\bP^1_m$ be a chain of $m$ copies of $\bP^1$. For $l=1,\ldots,m{-}1$, let
$q^1_l$ be the node at which $\bP^1_l$ and $\bP^1_{l+1}$ intersect. Let
$q^1_0\in \bP^1_1$ and $q^1_m\in \bP^1_m$ be smooth points.

A point in $\Mbar_{g,0}(\bP^1,\mu)$ is represented by a morphism
$$f\co(C,x_1,\ldots,x_h)\to (\bP^1[m], q^1_m)$$
where $C$ has at most nodal singularities, and $\bP^1[m]$ is obtained by
identifying $q^1\in\bP^1$ with $q^1_0\in\bP^1(m)$. In particular,
$\bP^1[0]=\bP^1$.  We call the original $\bP^1=\bP^1_0$ the root component
and $\bP^1_1,\ldots,\bP^1_m$ bubble components. We have $C=C_0\cup C_1
\cup \cdots \cup C_m$, where $C_l$ is the preimage of $\bP^1_l$. Let
$f_l\co C_l\to \bP^1_l$ be the restriction of $f$.  Then:
\begin{enumerate}
\item (\emph{degree})\qua $\deg f_l=d$, for $l=0,\ldots, h$.
\item (\emph{ramification})\qua  $f^{-1}(q^1_m)=\sum_{j=1}^h \mu_j x_j$.
\item (\emph{predeformability})\qua  The preimage of each node of the target
consists of nodes, at each of which two branches have the same contact
order. This is the predeformable condition: so that one can smooth both
the target and the domain to a morphism to $\bP^1$.
\item (\emph{stability})\qua  The automorphism group of $f$ is finite.
\end{enumerate}

Two morphisms satisfying (1)--(3) are equivalent if they have the same
target $\bP^1[m]$ for some nonnegative integer $m$ and they differ by
an isomorphism of the domain and an element of $\Aut(\bP^1[m],q^1_0,
q^1_m)\cong (\bC^*)^m$.  In particular, this defines the automorphism
group in (4). For fixed $g,\mu$, the stability condition (4) gives an
upper bound of the number $m$ of bubble components of the target.

By results in \cite{LiOne,LiTwo}, $\Mbar_{g,0}(\bP^1,\mu)$ is a proper,
separated  Deligne--Mumford stack with a perfect obstruction theory
of virtual dimension $2g{-}2{+}d{+}h$.  Roughly speaking, this means
that it is a compact, Hausdorff singular orbifold with a ``virtual
tangent bundle'' of rank $2g{-}2{+}d{+}h$. For later convenience, we will
also consider the disconnected version $\Mbar_\chi^\bu(\bP^1,\mu)$,
where the domain $C$ is allowed to be disconnected with
$2\bigl(h^0(\cO_C){-}h^1(\cO_C)\bigr)=\chi$.  Note that when $C$ is smooth,
$2\bigl(h^0(\cO_C){-}h^1(\cO_C)\bigr)$ is also the Euler characteristic
of $C$ as a compact surface. $\Mbar_\chi^\bu(\bP^1,\mu)$ is a proper,
separated Deligne--Mumford stack with a perfect obstruction theory of
virtual dimension $-\chi+d+h$.

Similarly, we may specify ramification types $\nu,\mu$ over two
points $q^0,q^1\in\bP^1$ and define the corresponding moduli
spaces
$\Mbar_{g,0}(\bP^1,\nu,\mu)$ and $\Mbar^\bu_\chi(\bP^1,\nu,\mu)$
of relative stable maps. We will also consider the quotient
$$\Mbar_\chi^\bu(\bP^1,\nu,\mu)//\bC^*\equiv
\bigl(\Mbar_\chi^\bu(\bP^1,\nu,\mu)\setminus
\Mbar_\chi^\bu(\bP^1,\nu,\mu)^{\bC^*}\bigr)/\bC^*$$
by the automorphism group $\bC^*$ of the target $(\bP^1,q^0,q^1)$.
\subsection{Double Hurwitz numbers as relative Gromov--Witten invariants}
The moduli space $\Mbar_\chi^\bu(\bP^1,\nu,\mu)$ parametrizes morphisms with targets 
of the form $\bP^1[m_0,m_1]$, where $\bP^1[m_0,m_1]$ is obtained by attaching $\bP^1(m_0)$ and 
$\bP^1(m_1)$ to $\bP^1$ at  $q^0$ and $q^1$ respectively. The distinguished points
on $\bP^1[m_0,m_1]$ are $q^0_{m_0}$ and $q^1_{m_1}$. Let $\pi_{m_0,m_1}\co\bP^1[m_0,m_1]\to \bP^1$ 
be the contraction to the root component.

There is a branch morphism
$$\Br\co\Mbar_\chi^\bu(\bP^1,\nu,\mu) \to \Sym^r\bP^1 \cong \bP^r$$
sending $\bigl[f\co C \to \bP^1[m_0,m_1] \bigr]$ to 
$$\mathrm{div}\bigl(\tf\,\bigr)- (d-\ell(\nu) ) q^0 -(d-\ell(\mu)) q^1$$
where $\mathrm{div}\bigl(\tf\,\bigr)$ is the branch divisor 
of $\tf=\pi_{m_0,m_1}\circ f\co C\to \bP^1$, and 
$$r=-\chi+\ell(\nu)+\ell(\mu)$$
is the virtual dimension of $\Mbar_\chi^\bu(\bP^1,\nu,\mu)$.

The structures on $\Mbar_\chi^\bu(\bP^1,\nu,\mu)$ allow one to
construct a \emph{virtual fundamental class}
$$\bigl[ \Mbar_\chi^\bu(\bP^1,\nu,\mu)\bigr]^\vir\in
  H_{2r}(\Mbar_\chi^\bu(\bP^1,\nu,\mu);\bQ)$$
which plays the role of the fundamental class of a compact oriented manifold.
Let $H\in H^2(\bP^r;\bZ)$ be the hyperplane class. The 
disconnected double Hurwitz number is equal to
$$
H^\bu_{\chi,\nu,\mu}=\frac{1}{|\Aut(\nu)|| \Aut(\mu) |}
\int_{[\Mbar_\chi^\bu(\bP^1,\nu,\mu)]^\vir}\Br^* H^r
$$
which is equal to
$$
\frac{r!}{|\Aut(\nu)|| \Aut(\mu) |}
\int_{[\Mbar_\chi^\bu(\bP^1,\nu,\mu)//\bC^*]^\vir}\psi_0^{r-1}
$$
by virtual localization, where $\psi_0$
is the \emph{target $\psi$--class}, the first Chern class
of the line bundle $\bL_0$ over $\Mbar_g(\bP^1,\nu,\mu)$
whose fiber at
$$
\bigl[f\co C\to \bP^1[m_0,m_1] \bigr]
$$ 
is the cotangent line $T^*_{q^0_{m_0}} \bP^1[m_0,m_1]$.

\subsection{Obstruction bundle}

We will define a vector bundle
$V^\bu_{\chi,\mu}$ over the moduli space $\Mbar^\bu_\chi(\bP^1,\mu)$.

For $f\co(C,x_1,\ldots, x_h)\to \bP^1[m]$, let $\tf$ be the composition $\pi_m\circ f\co C\to
\bP^1$, where $\pi_m\co\bP^1[m]\to \bP^1$ is the contraction to the 
root component. The fiber of $V^\bu_{\chi,\mu}$ at 
$$\bigl[ f\co(C,x_1,\ldots,x_h)\to \bP^1[m] \bigr]\in
\Mbar^\bu_\chi(\bP^1,\mu)$$
is 
\begin{equation}\label{eqn-Hone}
H^1(C,\tf^*\cO_{\bP^1}(-1))\oplus H^1(C,\cO_C(-x_1-\cdots-x_h)).
\end{equation}
Note that
$$H^0(C,\tf^*\cO_{\bP^1}(-1))= H^0(C,\cO_C(-x_1-\cdots-x_h))=0,$$
so \eqref{eqn-Hone} forms a vector bundle over $\Mbar^\bu_\chi(\bP^1,\mu)$
which, by Riemann--Roch, has rank $-\chi{+}d{+}h$, which is equal to
the virtual dimension of $\Mbar^\bu_\chi(\bP^1,\mu)$.  So
$$K^\bu_{\chi,\mu}=\frac{1}{|\Aut(\mu)|}
  \int_{\bigl[\Mbar_\chi^\bu(\bP^1,\mu)\bigr]^\vir}
  e\bigl(V^\bu_{\chi,\mu}\bigr)$$
is a topological invariant, where $e\bigl(V^\bu_{\chi,\mu}\bigr)$ is the Euler
class (top Chern class) of $V^\bu_{\chi,\mu}$.

The generating function in \fullref{thm-KGI} is given by
\begin{equation}\label{eqn-K}
K^\bu_\mu(\lam)=\sqrt{-1}^{\,h-d}\sum_{\chi}\lam^{-\chi+h}K^\bu_{\chi,\mu}.
\end{equation}

\subsection{Virtual localization}

Let $\bC^*$ act on $\bP^1$ by
$t\cdot [X,Y]=[tX,Y]$ for $t\in\bC^*$ and $[X,Y]\in\bP^1$. The two fixed
points are $q^0=[0,1]$ and $q^1=[1,0]$.  This induces a $\bC^*$--action on
$\Mbar^\bu_\chi(\bP^1,\mu)$. We would like to lift the $\bC^*$--action
to the vector bundle $V^\bu_{\chi,\mu}$ so that we can calculate
the integral of the equivariant Euler class $e_{\bC^*}(V^\bu_{\chi,\mu})$
by
virtual localization. It suffices to lift the $\bC^*$--action
on $\bP^1$ to the line bundles $\cO_{\bP^1}(-1)$ and $\cO_{\bP^1}$,
which will induce
actions on the cohomology groups \eqref{eqn-Hone}. We only need to specify
the weights of $\bC^*$--actions on the fibers of these line bundles
over $q^0$ and $q^1$:
$$
\begin{array}{l|ll}
     &   q^0 & q^1\\ \hline
\cO_{\bP^1}(-1) &  t^{-\tau-1} & t^{-\tau}\\
\cO_{\bP^1} &  t^{\tau} & t^{\tau}
\end{array}
$$
where $\tau\in\bZ$. We have
\begin{eqnarray*}
K_{\chi,\mu}^\bu
&=&\frac{1}{|\Aut(\mu)|}\int_{\bigl[\Mbar_\chi^\bu(\bP^1,\mu)\bigr]^\vir }
  e\bigl(V^\bu_{\chi,\mu}\bigr)\\
&=&\frac{1}{|\Aut(\mu)|}\int_{\bigl[\Mbar_\chi^\bu(\bP^1,\mu)\bigr]^\vir}
  e_{\bC^*}\bigl(V^\bu_{\chi,\mu}\bigr)\\
&=&\frac{1}{|\Aut(\mu)|}\sum_F\int_{[F]^\vir}
  \frac{e_{\bC^*}\bigl(V^\bu_{\chi,\mu}\bigr)}{e_{\bC^*}\bigl(N^\vir_F\bigr)}.
\end{eqnarray*}
If $\Mbar_\chi^\bu(\bP^1,\mu)$ is a compact complex manifold and each
fixed locus $F$ is a compact complex submanifold, then $[F]^\vir$ is 
the fundamental class of $F$, $N^\vir_F$ is the normal
bundle of $F$ in $\Mbar_\chi^\bu(\bP^1,\mu)$, and the last equality
is the Atiyah--Bott localization formula. Here we need to apply
virtual localization formula proved by Graber and Pandharipande in \cite{Gr-Pa},
where fundamental classes and normal bundles are replaced by
\emph{virtual} fundamental classes and \emph{virtual} normal bundles, respectively. 
The equivariant Euler class $e_{\bC^*}\bigl(V^\bu_{\chi,\mu}\bigr)$ and the
contribution from each fixed locus $F$ depend on $\tau$.

Let $f\co (C,x_1,\ldots,x_h)\to \bP^1[m]$ be a point in 
the fixed points set $\Mbar_\chi^\bu(\bP^1,\mu)^{\bC^*}$, let
$\tf=\pi_m\circ f\co C\to \bP^1$, and let $C^i=\tf^{-1}(q^i)$, for
$i=0,1$. Then the complement of $C^0\cup C^1$ in $C$ is a disjoint union
of twice-punctured spheres $L_1,\ldots, L_k$, and $f|_{L_j} \co L_j\to
\bP^1\setminus\{q^0,q^1\}$ is an honest covering map of some degree
$\nu_j$. This gives a partition $\nu=(\nu_1,\ldots,\nu_k)$ of $d=|\mu|$.
The restriction $f|_{C^0}$ is a constant map to $q^0$,  and $f|_{C^1}\co
C^1\to \bP^1(m)$ represents a point in
$\Mbar_{\chi^1}^\bu(\bP^1,\nu,\mu)//\bC^*$, so
the fixed locus is a finite quotient of 
\begin{equation}\label{eqn-fixed}
\Mbar_{\chi^0,\ell(\nu)}^\bu\times
\bigl(\Mbar_{\chi^1}^\bu(\bP^1,\nu,\mu)//\bC^*\bigr)
\end{equation}
where $\chi^0{+}\chi^1{-}2\ell(\nu)=\chi$, and
$\Mbar^\bu_{\chi,n}$ is the moduli stack of possibly disconnected stable
curves $C$ with $n$ marked points and $2(h^0(\cO_C){-}h^1(\cO_C))=\chi$.
The contribution from the above fixed locus is of the form
$$\int_{\Mbar^\bu_{\chi^0,\ell(\nu)}}\bigl(\cdots \bigr)
  \cdot A(\tau)\cdot
  \int_{\bigl[\Mbar^\bu_{\chi^1}(\bP^1,\nu,\mu)//\bC^*\bigr]^\vir}
  \bigl(\cdots\bigr)$$
where $A(\tau)$ is some combinatorial factor.  Calculations in \cite{LLZOne,LLZTwo} show that the integral over $\Mbar^\bu_{\chi^0,\ell(\nu)}$ is
exactly the one-partition Hodge integral $G^\bu_{\chi^0,\nu}(\tau)$
and the integrand of the integral over
$\Mbar^\bu_{\chi^1}(\bP^1,\nu,\mu)//\bC^*$ is a power of the target $\psi$
class. We have
\begin{equation}\label{eqn-KGH}
K^\bu_{\chi,\mu}=\sum_{\chi^0+\chi^1-2\ell(\nu)=\chi} 
G^\bu_{\chi^0,\nu}(\tau) \cdot A_{\chi^0,\nu,\chi^1}(\tau) \cdot H^\bu_{\chi^1,\nu,\mu}
\end{equation}
where $A_{\chi^0,\nu,\chi^1}(\tau)$ is some $\tau$--dependent
combinatorial factor. The expression \eqref{eqn-KGH} can be neatly
packaged in terms of generating functions:
$$K^\bu_{\chi}(\lam)=\sum_{|\nu|=d}
  G^\bu_{\nu}(\lam;\tau)z_\nu\Phi^\bu_{\nu,\mu}\bigl(-\sqrt{-1}\lam\tau\bigr).$$
This completes the proof of \fullref{thm-KGI}.

\section{Generalization to the two-partition case}
\label{sec-II}

\subsection{Gromov--Witten invariants of local toric Calabi--Yau threefolds}

Let $S$ be a Fano surface, and let $X$ be the total space
of the canonical line bundle $K_S$ of $S$ (for example, $\cO(-3)\to \bP^2$).
Then $X$ is a noncompact Calabi--Yau threefold. We call such
a noncompact Calabi--Yau threefold a \emph{local Calabi--Yau threefold}.
The image of any nonconstant morphism from a complex algebraic curve to $X$ is
contained in $S$, so for any nontrivial $d\in H_2(S;\bZ)\cong
H_2(X;\bZ)$ we have
$$\Mbar_{g,0}(X,d)=\Mbar_{g,0}(S,d)$$
as Deligne--Mumford stacks. However, they have different perfect
obstruction theories: the virtual dimension of the perfect obstruction theory
of $\Mbar_{g,0}(X,d)$ is zero, while that of $\Mbar_{g,0}(S,d)$ is
$$g-1-\int_{d} c_1(K_S).$$
The genus--$g$, degree--$d$ (for $d\neq0$) Gromov--Witten invariant of $X$
is defined by
$$N^X_{g,d}=\int_{\bigl[\Mbar_{g,0}(X,d)\bigr]^\vir} 1 \in\bQ.$$
Let $V_{g,d}$ be the vector bundle over $\Mbar_{g,0}(S,d)$ whose fiber
over a point represented by $f\co C\to S$ is given by
\begin{equation}\label{eqn-KS}
H^1(C,f^* K_S).
\end{equation}
Note that $H^0(C,f^* K_S)=0$, so \eqref{eqn-KS} forms a vector bundle
over $\Mbar_{g,0}(S,d)$. By Riemann--Roch, its rank is $g{-}1{-}\!\int_d
c_1(K_S)$ which is equal to the virtual dimension of $\Mbar_{g,0}(S,d)$. 
We have
\begin{equation}\label{eqn-NX}
N^X_{g,d}=\int_{ [\Mbar_{g,0}(X,d)]^\vir}1=\int_{ [\Mbar_{g,0}(S,d)]^\vir} e(V_{g,d}).
\end{equation}

When $X$ is a local \emph{toric} Calabi--Yau threefold, that is, the total
space of the canonical line bundle of a \emph{toric} Fano surface,
the integral in \eqref{eqn-NX} can be reduced to \emph{two-partition
Hodge integrals} by virtual localization.

\subsection{Two-partition Hodge integrals}
\label{sec-hodgeII}
Given $\mu^+,\mu^-$ partitions, let $\ell^\pm=\ell(\mu^\pm)$  
and define
\begin{eqnarray*}
G_{g,\mu^+,\mu^-}(\tau)&=&
  \frac{-\sqrt{-1}^{\,\ell^++\ell^-}}{|\Aut(\mu^+)||\Aut(\mu^-)|}
  (\tau(\tau+1))^{\ell^++\ell^--1}\\
&& \cdot \prod_{i=1}^{\ell^+}
   \frac{\prod_{a=1}^{\mu_i^+ -1}(\mu_i^+\tau+a)}{(\mu_i^+-1)!}
        \prod_{j=1}^{\ell^-}
   \frac{\prod_{a=1}^{\mu_j^- -1}\bigl(\frac{\mu_j^-}{\tau} + a\bigr)}{(\mu_j^-
     -1)!}\\
&& \cdot\int_{\Mbar_{g,\ell^+ +\ell^-}}
\frac{\Lambda^\vee_g(1)\Lambda^\vee_g(\tau)\Lambda^\vee_g(-\tau-1)}
  {\prod_{i=1}^{\ell^+}\bigl(1-\mu_i^+\psi_i\bigr)
  \prod_{j=1}^{\ell^-}\tau\bigl(\tau-\mu_j^-\psi_{\ell^++j}\bigr)}
\end{eqnarray*}
In particular,
\begin{equation}\label{eqn-twoone}
\sqrt{-1}^{|\mu|}G_{g,\mu,\emptyset}(\tau)=G_{g,\mu}(\tau)
\end{equation}
where $G_{g,\mu}(\tau)$ is the one-partition Hodge integral defined in
\fullref{sec-hodgeI}.

\subsection{Formula of two-partition Hodge integrals}
To state the formula of two-partition Hodge integrals, we introduce
some generating functions.

We first define generating function of two-partition Hodge integrals. 
Introduce variables 
$$p^+=(p_1^+,p_2^+,\ldots) \quad\text{and}\quad p^-=(p_1^-,p_2^-,\ldots).$$
Given a partition $\mu$, define
$$p^\pm_\mu=p^\pm_1\ldots p^\pm_{\ell(\mu)}$$
Define generating functions
\begin{eqnarray*}
G_{\mu^+,\mu^-}(\lambda;\tau)&=&
\sum_{g=0}^\infty \lambda^{2g-2+\ell^++\ell^-}G_{g,\mu^+,\mu^-}(\tau)\\
G(\lambda;p^+,p^-;\tau)&=&\sum_{(\mu^+,\mu^-)\neq(\emptyset,\emptyset)}
G_{\mu^+,\mu^-}(\lambda;\tau)p^+_{\mu^+}p^-_{\mu^-}\\
G^\bullet(\lambda;p^+,p^-;\tau)&=&
  \exp\bigl(G\bigl(\lambda;p^+,p^-;\tau\bigr)\bigr)\\
&=& \sum_{\mu^+,\mu^-}G^\bu_{\mu^+,\mu^-}(\lam;\tau)p^+_{\mu^+}p^-_{\mu^-}
\end{eqnarray*}
$G^\bu_{\mu^+,\mu^-}(\lam;\tau)$ is the disconnected version of 
$G_{\mu^+,\mu^-}(\lam;\tau)$. By \eqref{eqn-twoone}, 
\begin{equation}\label{eqn-GIII}
\sqrt{-1}^{|\mu|}G^\bu_{\mu,\emptyset}(\lam;\tau)=G^\bu_\mu(\lam;\tau).
\end{equation}
We next define generating functions of symmetric group
representations. Let 
$$s_\mu(x_1,x_2,\ldots)$$
be \emph{Schur functions} (see \cite{Ma} for definitions). Recall that
$$s_\nu s_\rho=\sum_\mu c^\mu_{\nu\rho}s_\mu \quad\text{and}\quad
s_{\mu/\rho}=\sum_\rho c^{\mu}_{\nu\rho}s_\rho$$
where $c^\mu_{\nu\rho}$ are known as
\emph{Littlewood--Richardson coefficients} and $s_{\mu/\rho}$ are known as
\emph{skew Schur functions}. Let $q=e^{\sqrt{-1}\lam}$ and
write
$$q^{-\rho}=\bigl(q^\frac{1}{2}, q^\frac{3}{2},\ldots\bigr)$$
Introduce
$$\cW_{\mu,\nu}(q)=(-1)^{|\mu|+|\nu|} q^{\frac{\kappa_\mu+\kappa_\nu}{2}}
\sum_\eta s_{\mu/\eta}(q^{-\rho})s_{\nu/\eta}(q^{-\rho})$$
which is related to the HOMFLY polynomial of the Hopf link.
In particular,
$$
\cW_{\mu,\emptyset}(q)=(-1)^{|\mu|}q^{\kappa_\mu/2}s_\mu(q^{-\rho})=s_\mu(q^\rho)=\cW_\mu(q)
$$
where $\cW_\mu(q)$ is defined by \eqref{eqn-WI} (see \cite{ZhFour} for
details). Define
\begin{equation}\label{eqn-Rtwo}
R^\bu_{\mu^+,\mu^-}(\lambda;\tau)
=\sum_{|\nu^\pm|=|\mu^\pm|}\frac{\chi_{\nu^+}(C_{\mu^+})}{z_{\mu^+}}
     \frac{\chi_{\nu^-}(C_{\mu^-})}{z_{\mu^-}} e^{\sqrt{-1}(\kappa_{\nu^+}\tau +\kappa_{\nu^-}\tau^{-1})\lambda/2}
\mathcal{W}_{\nu^+,\nu^-}(q).
\end{equation}
In particular,
\begin{equation}\label{eqn-RIII}
\sqrt{-1}^{|\mu|}R^\bu_{\mu,\emptyset}(\lam;\tau)=
\sqrt{-1}^{|\mu|}\sum_\nu \frac{\chi_\nu(C_\mu)}{z_\mu}e^{\sqrt{-1}\kappa_\nu\tau\lam/2}\cW_\nu(q)= R^\bu_\mu(\lam;\tau)
\end{equation}
where $R^\bu_\mu(\lam;\tau)$  is defined by \eqref{eqn-RV}.
\begin{conj}[Zhou \cite{ZhTwo}]\label{conj-GRtwo}
\begin{equation}\label{eqn-GRtwo}
 G^\bu_{\mu^+,\mu^-}(\lam;\tau)=R^\bu_{\mu^+,\mu^-}(\lam;\tau).
\end{equation}
\end{conj}
By \eqref{eqn-GIII} and \eqref{eqn-RIII},
\fullref{conj-GRtwo} reduces to the Mari\~{n}o--Vafa formula when 
$\mu^-=\emptyset$. 

\subsection{Application}

Let $X$ be a local toric Calabi--Yau threefold.  The \emph{Gromov--Witten
potential} of $X$ is defined by
$$F^X(\lam,t)=\sum_{g=0}^\infty
  \sum_{\tiny\begin{array}{c}d\in H_2(X,\bZ)\\d\neq 0\end{array}}
  \lam^{2g-2}N^X_{g,d} e^{-d\cdot t}$$
where $t=(t_1,t_2,\ldots)$ are coordinates on $H^{1,1}(X)$. Note that
$H^2(X,\bC)=H^{1,1}(X)$
because $X$ is toric. The \emph{partition function} is defined by
$$Z^X(\lam,t)=\exp\bigl(F^X(\lam,t)\bigr)$$
By virtual localization, $Z^X(\lam,t)$ can be expressed in terms of
$G^\bu_{\mu^+,\mu^-}(\lam;\tau)$, so \eqref{eqn-GRtwo}
gives an explicit formula of $Z^X(\lam,t)$ for any local toric Calabi--Yau
threefold $X$ in terms of $\cW_{\mu\nu}(q)$.   

For example, let $X$ be the total space of $\cO_{\bP^1}(-3)\to \bP^2$. Then
\begin{align*}
Z^X(\lam,\!t)=\!\!\smash{\sum_{\mu^i\neq \emptyset}
  \bigl((-1)^{\sum_{i=1}^3|\mu^i|}\bigr)
  \bigl(e^{-\sum_{i=1}^3|\mu^i|t} \bigr)}&
  \bigl(q^{\sum_{i=1}^3\kappa_{\nu^i}/2} \bigr)\\
  &\cW_{\mu^1\mu^2}(q)\cW_{\mu^2\mu^3}(q)\cW_{\mu^3\mu^1}(q).
\end{align*}
The algorithm for computing $Z^X(\lam,t)$ for local toric
Calabi--Yau threefolds in terms of $\cW_{\mu\nu}(q)$ was described by
Aganagic--Mari\~{n}o--Vafa \cite{Ag-Ma-Va}, and an explicit formula
was given by Iqbal in \cite{Iq}.  Motivated by the Mari\~{n}o--Vafa
and Iqbal formul\ae, Zhou conjectured a formula \eqref{eqn-GRtwo}
for two-partition Hodge integrals and used it, together with virtual
localization, to give a mathematical derivation of Iqbal's formula. The
formula for two-partition Hodge integrals was proved in \cite{LLZFour}.

\subsection{Outline of proof}
In this subsection,  we outline the proof of the formula of
two-partition Hodge integrals given in \cite{LLZFour}. 
The strategy is similar to the second approach to the Mari\~{n}o--Vafa
formula described in \fullref{sec-dHurwitz}. 

It follows from \eqref{eqn-Rtwo} that
\begin{equation}\label{eqn-RRtwo}
\begin{aligned}
R&^\bu_{\mu^+,\mu^-}(\lam;\!\tau) = \\
&\sum_{|\nu^\pm|=|\mu^\pm|}\!\!\!\!R^\bu_{\mu^+,\mu^-}(\lam;\!\tau_0)
   {\cdot}
z_{\nu^+}\Phi^\bu_{\nu^+,\mu^+}\bigl(\sqrt{-1}\lam(\tau{-}\tau_0)\bigr)
{\cdot} z_{\nu^-}\Phi^\bu_{\nu^-,\mu^-}
\bigl(\sqrt{-1}\lam\bigl({\textstyle\frac{1}{\tau}}{-}
  {\textstyle\frac{1}{\tau_0}} \bigr)\bigr)
\end{aligned}
\end{equation}
for any $\tau \in \bC^*$. Here we cannot specialize at $\tau=0$
because both sides of \eqref{eqn-GRtwo} have a pole at $\tau=0$.
At $\tau=-1$, the two partition integrals vanish
unless $\ell(\mu^+)+\ell(\mu^-)=1$ and we are left with
$b_g$ (the $\lambda_g$--integrals on $\Mbar_{g,1}$).  
\begin{theo}[Zhou \cite{ZhTwo}]
\begin{equation}\label{eqn-initialtwo}
G^\bu_{\mu^+,\mu^-}(\lambda;-1)=R^\bu_{\mu^+,\mu^-}(\lambda;-1)
\end{equation}
\end{theo}

The authors of \cite{LLZFour} defined 
generating functions $K^\bu_{\mu^+,\mu^-}(\lam)$ for
relative Gromov--Witten invariants of $\bP^1\times \bP^1$ blown up at a point, 
and used virtual localization to derive the following expression.
\begin{prop}[Liu--Liu--Zhou \cite{LLZFour}]\label{thm-KGtwo}
\begin{equation}
K^\bu_{\mu^+,\mu^-}(\lam)
=\!\!\!\!\sum_{|\nu^\pm|=|\mu^\pm|}\!\!\!\!G^\bu_{\mu^+,\mu^-}(\lam;\tau)
z_{\nu^+}\Phi^\bu_{\nu^+,\mu^+}\bigl(-\sqrt{-1}\lam\tau\bigr)
z_{\nu^-}\Phi^\bu_{\nu^-,\mu^-}\bigl(-\sqrt{-1}\lam/\tau\bigr)
\end{equation}
\end{prop}

This proposition implies the following convolution equation 
of $G^\bu_{\mu^+,\mu^-}(\lam;\tau)$:
\begin{theo}
\begin{equation}\label{eqn-GGtwo}
\begin{aligned}
&G^\bu_{\mu^+,\mu^-}(\lam;\!\tau) = \\
&\sum_{|\nu^\pm|=|\mu^\pm|}\!\!\!\!G^\bu_{\mu^+,\mu^-}(\lam;\!\tau_0)
{\cdot} z_{\nu^+}\Phi^\bu_{\nu^+,\mu^+}(\sqrt{-1}\lam\bigl(\tau{-}\tau_0)\bigr)
{\cdot} z_{\nu^-}\Phi^\bu_{\nu^-,\mu^-}
\bigl(\sqrt{-1}\lam\bigl({\textstyle\frac{1}{\tau}}{-}\textstyle{\frac{1}{\tau_0}}\bigr)\lam\bigr)
\end{aligned}
\end{equation}
for any $\tau_0\in \bC^*$.
\end{theo}
The formula \eqref{eqn-GRtwo} for two-partition Hodge integrals follows from
the convolution equation \eqref{eqn-RRtwo} for $R^\bu_{\mu^+,\mu^-}(\lam;\tau)$,
the convolution equation \eqref{eqn-GGtwo} for $G^\bu_{\mu^+,\mu^-}(\lam;\tau)$,
and the initial values \eqref{eqn-initialtwo}.

\bibliographystyle{gtart}
\bibliography{link}

\begin{thebibliography}{}
\providecommand\bibmarginpar{\leavevmode\marginpar}
\def\urlstyle#1{{\tt #1}}

\bibitem{AKMV}
\textbf{M Aganagic}, \textbf{A Klemm}, \textbf{M Mari{\~n}o}, \textbf{C Vafa},
  \emph{The topological vertex} \xox{arXiv}{hep-th/0305132}

\bibitem{Ag-Ma-Va}
\textbf{M Aganagic}, \textbf{M Mari{\~n}o}, \textbf{C Vafa},
  \href{http://dx.doi.org/10.1007/s00220-004-1067-x} {\emph{All loop
  topological string amplitudes from {C}hern--{S}imons theory}}, Comm. Math.
  Phys. 247 (2004) 467--512 \xox{MR}{2063269}

\bibitem{Di-Fl}
\textbf{D-E Diaconescu}, \textbf{B Florea}, \emph{Localization and gluing of
  topological amplitutes} \xox{arXiv}{hep-th/0309143}

\bibitem{Di-Fl-Gr}
\textbf{D-E Diaconescu}, \textbf{B Florea}, \textbf{A Grassi}, \emph{Geometric
  transitions, del {P}ezzo surfaces and open string instantons}, Adv. Theor.
  Math. Phys. 6 (2002) 643--702 \xox{MR}{1969655}

\bibitem{ELSVOne}
\textbf{T Ekedahl}, \textbf{S Lando}, \textbf{M Shapiro}, \textbf{A
  Vainshtein}, \href{http://dx.doi.org/10.1016/S0764-4442(99)80435-2} {\emph{On
  {H}urwitz numbers and {H}odge integrals}}, C. R. Acad. Sci. Paris S\'er. I
  Math. 328 (1999) 1175--1180 \xox{MR}{1701381}

\bibitem{ELSVTwo}
\textbf{T Ekedahl}, \textbf{S Lando}, \textbf{M Shapiro}, \textbf{A
  Vainshtein}, \href{http://dx.doi.org/10.1007/s002220100164} {\emph{Hurwitz
  numbers and intersections on moduli spaces of curves}}, Invent. Math. 146
  (2001) 297--327 \xox{MR}{1864018}

\bibitem{Fa-PaOne}
\textbf{C Faber}, \textbf{R Pandharipande}, \emph{Hodge integrals and
  {G}romov--{W}itten theory}, Invent. Math. 139 (2000) 173--199
  \xox{MR}{1728879}

\bibitem{Fa-PaTwo}
\textbf{C Faber}, \textbf{R Pandharipande}, \emph{Hodge integrals, partition
  matrices, and the $\lambda_g$ conjecture}, Ann. of Math. $(2)$ 157 (2003)
  97--124 \xox{MR}{1954265}

\bibitem{Ge-Pa}
\textbf{E Getzler}, \textbf{R Pandharipande},
  \href{http://dx.doi.org/10.1016/S0550-3213(98)00517-3} {\emph{Virasoro
  constraints and the {C}hern classes of the {H}odge bundle}}, Nuclear Phys. B
  530 (1998) 701--714 \xox{MR}{1653492}

\bibitem{Gi}
\textbf{A\,B Givental}, \emph{Gromov--{W}itten invariants and quantization of
  quadratic {H}amiltonians}, Mosc. Math. J. 1 (2001) 551--568, 645
  \xox{MR}{1901075}

\bibitem{Go-Va}
\textbf{R Gopakumar}, \textbf{C Vafa}, \emph{On the gauge theory/geometry
  correspondence}, Adv. Theor. Math. Phys. 3 (1999) 1415--1443
  \xox{MR}{1796682}

\bibitem{Go-Ja-Va}
\textbf{I\,P Goulden}, \textbf{D\,M Jackson}, \textbf{A Vainshtein},
  \href{http://dx.doi.org/10.1007/PL00001274} {\emph{The number of ramified
  coverings of the sphere by the torus and surfaces of higher genera}}, Ann.
  Comb. 4 (2000) 27--46 \xox{MR}{1763948}

\bibitem{Gr-Pa}
\textbf{T Graber}, \textbf{R Pandharipande},
  \href{http://dx.doi.org/10.1007/s002220050293} {\emph{Localization of virtual
  classes}}, Invent. Math. 135 (1999) 487--518 \xox{MR}{1666787}

\bibitem{Gr-VaOne}
\textbf{T Graber}, \textbf{R Vakil},
  \href{http://dx.doi.org/10.1023/A:1021791611677} {\emph{Hodge integrals and
  {H}urwitz numbers via virtual localization}}, Compositio Math. 135 (2003)
  25--36 \xox{MR}{1955162}

\bibitem{Io-PaTwo}
\textbf{E-N Ionel}, \textbf{T\,H Parker}, \emph{The symplectic sum formula for
  Gromov--Witten invariants} \xox{arXiv}{math.SG/0010217}

\bibitem{Io-PaOne}
\textbf{E-N Ionel}, \textbf{T\,H Parker}, \emph{Relative {G}romov--{W}itten
  invariants}, Ann. of Math. $(2)$ 157 (2003) 45--96 \xox{MR}{1954264}

\bibitem{Iq}
\textbf{A Iqbal}, \emph{All genus topological amplitudes and 5--brane webs as
  Feynman diagrams} \xox{arXiv}{hep-th/0207114}

\bibitem{Ka-Li}
\textbf{S Katz}, \textbf{C-C\,M Liu},
  \href{http://dx.doi.org/10.2140/gtm.2006.8.1} {\emph{Enumerative geometry of
  stable maps with {L}agrangian boundary conditions and multiple covers of the
  disc}}, from: ``The interaction of finite-type and Gromov-Witten invariants
  (Banff 2003)'', Geom. Topol. Monogr. 8 (2006)  1--48

\bibitem{KoOne}
\textbf{M Kontsevich},
  \href{http://projecteuclid.org/getRecord?id=euclid.cmp/1104250524}
  {\emph{Intersection theory on the moduli space of curves and the matrix
  {A}iry function}}, Comm. Math. Phys. 147 (1992) 1--23 \xox{MR}{1171758}

\bibitem{La-Ma-Va}
\textbf{J\,M\,F Labastida}, \textbf{M Mari{\~n}o}, \textbf{C Vafa},
  \href{http://dx.doi.org/10.1088/1126-6708/2000/11/007} {\emph{Knots, links
  and branes at large $N$}}, J. High Energy Phys.  (2000) Paper 7, 42 pages
  \xox{MR}{1806596}

\bibitem{Li-Ru}
\textbf{A-M Li}, \textbf{Y Ruan},
  \href{http://dx.doi.org/10.1007/s002220100146} {\emph{Symplectic surgery and
  {G}romov--{W}itten invariants of {C}alabi--{Y}au 3--folds}}, Invent. Math.
  145 (2001) 151--218 \xox{MR}{1839289}

\bibitem{Li-Zh-Zh}
\textbf{A-M Li}, \textbf{G Zhao}, \textbf{Q Zheng},
  \href{http://dx.doi.org/10.1007/s002200000254} {\emph{The number of ramified
  covering of a {R}iemann surface by {R}iemann surface}}, Comm. Math. Phys. 213
  (2000) 685--696 \xox{MR}{1785434}

\bibitem{LiOne}
\textbf{J Li},
  \href{http://projecteuclid.org/getRecord?id=euclid.jdg/1090348132}
  {\emph{Stable morphisms to singular schemes and relative stable morphisms}},
  J. Differential Geom. 57 (2001) 509--578 \xox{MR}{1882667}

\bibitem{LiTwo}
\textbf{J Li},
  \href{http://projecteuclid.org/getRecord?id=euclid.jdg/1090351102} {\emph{A
  degeneration formula of {GW}-invariants}}, J. Differential Geom. 60 (2002)
  199--293 \xox{MR}{1938113}

\bibitem{LLLZ}
\textbf{J Li}, \textbf{C-C\,M Liu}, \textbf{K Liu}, \textbf{J Zhou}, \emph{A
  mathematical theory of the topological vertex} \xox{arXiv}{math.AG/0408426}

\bibitem{Li-So}
\textbf{J Li}, \textbf{Y\,S Song},
  \href{http://dx.doi.org/10.2140/gtm.2006.8.49} {\emph{Open string instantons
  and relative stable morphisms}}, from: ``The interaction of finite-type and
  Gromov--Witten invariants (Banff 2003)'', Geom. Topol. Monogr. 8 (2006)
  49--72

\bibitem{LLY}
\textbf{B\,H Lian}, \textbf{K Liu}, \textbf{S-T Yau}, \emph{Mirror principle
  III}, Asian J. Math. 3 (1999) 771--800 \xox{MR}{1797578}

\bibitem{LLZFour}
\textbf{C-C\,M Liu}, \textbf{K Liu}, \textbf{J Zhou}, \emph{A formula of
  two-partition Hodge integrals} \xox{arXiv}{math.AG/0310272}

\bibitem{LLZThree}
\textbf{C-C\,M Liu}, \textbf{K Liu}, \textbf{J Zhou}, \emph{Mari{\~n}o--Vafa
  formula and Hodge integral identities} \xox{arXiv}{math.AG/0308015}

\bibitem{LLZTwo}
\textbf{C-C\,M Liu}, \textbf{K Liu}, \textbf{J Zhou},
  \href{http://projecteuclid.org/getRecord?id=euclid.jdg/1090511689} {\emph{A
  proof of a conjecture of {M}ari{\~n}o--{V}afa on {H}odge integrals}}, J.
  Differential Geom. 65 (2003) 289--340 \xox{MR}{2058264}

\bibitem{LLZOne}
\textbf{C-C\,M Liu}, \textbf{K Liu}, \textbf{J Zhou}, \emph{On a proof of a
  conjecture of {M}ari{\~n}o--{V}afa on {H}odge integrals}, Math. Res. Lett. 11
  (2004) 259--272 \xox{MR}{2067471}

\bibitem{Ma}
\textbf{I\,G Macdonald}, \emph{Symmetric functions and {H}all polynomials},
  Oxford Mathematical Monographs, The Clarendon Press Oxford University Press,
  New York (1995) \xox{MR}{1354144}

\bibitem{MaTwo}
\textbf{M Mari{\~n}o}, \emph{Chern--Simons theory and topological strings}
  \xox{arXiv}{hep-th/0406005}

\bibitem{MaOne}
\textbf{M Mari{\~n}o}, \emph{Enumerative geometry and knot invariants}
  \xox{arXiv}{hep-th/0210145}

\bibitem{Ma-Va}
\textbf{M Mari{\~n}o}, \textbf{C Vafa}, \emph{Framed knots at large $N$}, from:
  ``Orbifolds in mathematics and physics (Madison, WI, 2001)'', Contemp. Math.
  310, Amer. Math. Soc., Providence, RI (2002)  185--204 \xox{MR}{1950947}

\bibitem{MNOPOne}
\textbf{D Maulik}, \textbf{N Nekrasov}, \textbf{A Okounkov}, \textbf{R
  Pandharipande}, \emph{Gromov--Witten theory and Donaldson--Thomas theory I}
  \xox{arXiv}{math.AG/0312059}

\bibitem{MNOPTwo}
\textbf{D Maulik}, \textbf{N Nekrasov}, \textbf{A Okounkov}, \textbf{R
  Pandharipande}, \emph{Gromov--Witten theory and Donaldson--Thomas theory II}
  \xox{arXiv}{math.AG/0406092}

\bibitem{Mi}
\textbf{M Mirzakhani}, \emph{Weil--Petersson volumes and the Witten--Kontsevich
  formula} (2003)

\bibitem{Mu}
\textbf{D Mumford}, \emph{Towards an enumerative geometry of the moduli space
  of curves}, from: ``Arithmetic and geometry, Vol. II'', Progr. Math. 36,
  Birkh\"auser, Boston (1983)  271--328 \xox{MR}{717614}

\bibitem{Ok-PaOne}
\textbf{A Okounkov}, \textbf{R Pandharipande}, \emph{Gromov--Witten theory,
  Hurwitz numbers, and matrix models I} \xox{arXiv}{math.AG/0101147}

\bibitem{Ok-PaFour}
\textbf{A Okounkov}, \textbf{R Pandharipande}, \emph{Virasoro constraints for
  target curves} \xox{arXiv}{math.AG/0308097}

\bibitem{Ok-PaThree}
\textbf{A Okounkov}, \textbf{R Pandharipande},
  \href{http://dx.doi.org/10.2140/gt.2004.8.675} {\emph{Hodge integrals and
  invariants of the unknot}}, Geom. Topol. 8 (2004) 675--699 \xox{MR}{2057777}

\bibitem{Oo-Va}
\textbf{H Ooguri}, \textbf{C Vafa},
  \href{http://dx.doi.org/10.1016/S0550-3213(00)00118-8} {\emph{Knot invariants
  and topological strings}}, Nuclear Phys. B 577 (2000) 419--438
  \xox{MR}{1765411}

\bibitem{Ta}
\textbf{C\,H Taubes}, \emph{Lagrangians for the {G}opakumar--{V}afa
  conjecture}, Adv. Theor. Math. Phys. 5 (2001) 139--163 \xox{MR}{1894340}

\bibitem{WiOne}
\textbf{E Witten}, \emph{Two-dimensional gravity and intersection theory on
  moduli space}, from: ``Surveys in differential geometry (Cambridge, MA,
  1990)'', Lehigh Univ., Bethlehem, PA (1991)  243--310 \xox{MR}{1144529}

\bibitem{WiTwo}
\textbf{E Witten}, \emph{Chern--{S}imons gauge theory as a string theory},
  from: ``The Floer memorial volume'', Progr. Math. 133, Birkh\"auser, Basel
  (1995)  637--678 \xox{MR}{1362846}

\bibitem{ZhTwo}
\textbf{J Zhou}, \emph{A conjecture on Hodge integrals}
  \xox{arXiv}{math.AG/0310282}

\bibitem{ZhFour}
\textbf{J Zhou}, \emph{Curve counting and instanton counting}
  \xox{arXiv}{math.AG/0311237}

\bibitem{ZhOne}
\textbf{J Zhou}, \emph{Hodge integrals, Hurwitz numbers, and symmetric groups}
  \xox{arXiv}{math.AG/0308024}

\bibitem{ZhThree}
\textbf{J Zhou}, \emph{Localizations on moduli spaces and free field
  realizations of Feynman rules} \xox{arXiv}{math.AG/0310283}

\end{thebibliography}

\end{document}